\documentclass{amsart}
\usepackage{amsmath,amsfonts,amsthm,mathrsfs,amscd}
\usepackage[numeric,initials]{amsrefs}
\usepackage{amssymb}
\usepackage[PostScript=dvips,dpi=1270,heads=LaTeX]{diagrams} 
\diagramstyle{small,thin,labelstyle=\scriptstyle}
\newif\ifcras\crasfalse

\newtheorem{thm}{Theorem}

\newtheorem{e-proposition}[thm]{Proposition}
\newtheorem{prop}[thm]{Proposition}
\newtheorem*{lem*}{Lemma}
\newtheorem{cor}[thm]{Corollary}
\newtheorem*{thm*}{Theorem}
\newtheorem*{prop*}{Proposition}
\theoremstyle{definition}

\newcommand{\QQ}{\mathbb{Q}}
\newcommand{\ZZ}{\mathbb{Z}}
\newcommand{\CC}{\mathbb{C}}

\newcommand{\shO}{\mathscr{O}}
\newcommand{\define}[1]{\emph{#1}}
\newcommand{\Sb}{\bar{S}}
\newcommand{\shV}{\mathscr{V}}
\newcommand{\shVZ}{\shV_{\ZZ}}
\newcommand{\shVO}{\shV_{\shO}}
\newcommand{\VC}{V_{\CC}}
\newcommand{\VZ}{V_{\ZZ}}
\newcommand{\shH}{\mathscr{H}}
\DeclareMathOperator{\Hdg}{Hdg}

\DeclareMathOperator{\End}{End}

\newcommand{\Gr}{\mathit{Gr}}

\DeclareMathOperator{\Ext}{Ext}
\newcommand{\MHS}{\mathrm{MHS}}
\newcommand{\VMHS}{\mathrm{VMHS}}
\DeclareMathOperator{\NF}{NF}
\newcommand{\Yb}{\bar{Y}}
\DeclareMathOperator{\CH}{CH}
\newcommand{\Hmot}{H_{\mathcal{M}}}
\newcommand{\Hdel}{H_{\mathcal{D}}}
\DeclareMathOperator{\MHM}{MHM}

\newcommand{\pil}{\pi_{\ast}}

\newcommand{\tensor}{\otimes}
\DeclareMathOperator{\SL}{SL}
\renewcommand{\define}[1]{\emph{#1}}

\begin{document}
\title[The locus of Hodge classes]{The locus of Hodge classes in an admissible
variation of mixed Hodge structure}
\author[P.~Brosnan]{Patrick Brosnan}
\address{Department of Mathematics \\
The University of British Columbia \\
1984 Mathematics Road \\
Vancouver, B.C., Canada V6T 1Z2}
\email{brosnan@math.ubc.ca}
\author[G.~Pearlstein]{Greg Pearlstein}
\address{Department of Mathematics \\
Michigan State University \\
East Lansing, MI 48824}
\email{gpearl@math.msu.edu}
\author[C.~Schnell]{Christian Schnell}
\address{Department of Mathematics, Statistics \& Computer Science \\
University of Illinois at Chicago \\
851 South Morgan Street \\
Chicago, IL 60607}
\email{cschnell@math.uic.edu}
\begin{abstract}
We generalize the theorem of E.~Cattani, P.~Deligne, and A.~Kaplan to admissible
variations of mixed Hodge structure.
\end{abstract}
\maketitle

\section{Introduction} 

The purpose of this note is to prove the following generalization of the famous
theorem of Cattani, Deligne, and Kaplan \cite{CDK}.

\begin{thm} \label{thm:main}
Let $S$ be a Zariski-open subset of a complex manifold $\Sb$, and let $\shV$ be a
variation of mixed Hodge structure on $S$. Suppose that $\shV$ is defined over $\ZZ$,
graded polarized, and admissible with respect to $\Sb$. Let $\Hdg(\shV)$ denote the
locus of Hodge classes in $\shV$.  Then each component of $\Hdg(\shV)$ extends to an
analytic space, finite and proper over $\Sb$.
\end{thm}

As in the original paper, Chow's theorem implies that the locus of Hodge classes
consists of algebraic varieties if $S$ is algebraic.

\begin{cor}
In the situation of Theorem~\ref{thm:main}, suppose that $S$ is quasi-projective.
Then each component of $\Hdg(\shV)$ is a quasi-projective algebraic variety.
\end{cor}

We remind the reader of a few basic definitions. Given a mixed Hodge structure $V$
defined over $\ZZ$, a \define{Hodge class} in $V$ is an element of $\VZ \cap W_0 \VC
\cap F^0 \VC$, or equivalently, a morphism of mixed Hodge structures $\ZZ(0) \to V$. Given a
variation of mixed Hodge structure $\shV$ on a complex manifold $S$, let $\shVZ$ denote the
underlying integral local system. Its \'etal\'e space $T(\shVZ)$ is a covering space
of $S$ with countably many connected components; it naturally embeds into the
holomorphic vector bundle $E(\shVO)$. The \define{locus of Hodge classes} in $\shV$ can then be
described as the intersection
\[
	\Hdg(\shV) = T(\shVZ) \cap E \bigl( F^0 \shVO \bigr) \cap E \bigl( W_0 \shVO \bigr).
\]

We deduce Theorem~\ref{thm:main} from the original result by Cattani, Deligne, and
Kaplan with the help of the following difficult theorem; it is the main result of
\cite{BP}, and can also be proved by the methods of \cite{Schnell-Neron}. (A similar
result has also been announced by Kato, Nakayama, and Usui in \cite{KNU-mod}.) Either
proof relies on the $\SL(2)$-orbit theorem of Kato, Nakayama, and Usui \cite{KNU}. 

\begin{thm} \label{thm:BP}
Let $\nu$ be an admissible higher normal function on $S$, that is, an admissible
extension of $\ZZ(0)$ by a polarized variation of Hodge structure of
negative weight.
Let $Z(\nu)=\{s\in S: \nu(s)=0\}$ denote the zero locus of
$\nu$. (C.f. the discussion at the beginning of Section~\ref{s.LocusOfHodge}.)
Then the closure of $Z(\nu)$ in $\Sb$ is an analytic subset.
\end{thm}

Note that this result includes the case of classical normal functions (where the
Hodge structure has weight $-1$).  Theorem~\ref{thm:BP} in itself is most interesting
when $S$ is a quasi-projective complex manifold; we may then take $\Sb$ to be any
smooth projective compactification, since the notion of admissibility is independent
of the particular choice.

\begin{cor} \label{cor:HNF}
Suppose that $\nu$ is an admissible higher normal function on $S$, that is, an extension of
$\ZZ(0)$ by a polarized variation of Hodge structure of negative weight. Then the
zero locus $Z(\nu)$ is an algebraic subset of $S$.
\end{cor}

One source for higher normal functions is through families of higher Chow
cycles. Let $\pi \colon X \to S$ be a family of complex projective manifolds with $S$
smooth. Then the regulator map from motivic cohomology $\Hmot^p \bigl( X, \ZZ(q)
\bigr) \simeq \CH^q(X, 2q-p)$ to Deligne cohomology $\Hdel^p \bigl( X, \ZZ(q) \bigr)$
induces a homomorphism
\[
	\CH^q(X, 2q-p) \tensor \QQ \to \bigoplus_{k \in \ZZ} 
		\Ext_{\MHM(S)}^{p-k} \bigl( \QQ(0), R^k \pil \QQ(q) \bigr),
\]
using the decomposition theorem; $\MHM(S)$ is the category of mixed Hodge
modules on $S$. In particular, a higher Chow cycle on $X$ determines an element in
$\Ext_{\MHM(S)}^1 \bigl( \QQ, R^{p-1} \pil \QQ(q) \bigr)$; some multiple is an
admissible higher normal function for the variation of Hodge structure $R^{p-1} \pil
\ZZ(q)$ of weight $p-2q-1 < 0$.

The same methods can be used to describe the locus of points $s \in S$ where $V_s$
splits over $\ZZ$ (we say that a mixed Hodge structure $V$ \define{splits} 
over $\ZZ$ if $V \simeq \bigoplus_m \Gr_m^W V$ in $\MHS$).

\begin{thm} \label{thm:split}
Let $\shV$ be an admissible variation of mixed Hodge structure on $S$. Then the set
of points $s \in S$ where the mixed Hodge structure $V_s$ splits over $\ZZ$ is an
algebraic subset of $S$.
\end{thm}

Since $V_s$ splits over $\ZZ$ iff there is a Hodge class in $\End(V_s)$ that induces
a splitting of the underlying integral lattice, this result may also be viewed as a
special case of Theorem~\ref{thm:main}.

\section{Admissibility}

Let $\shV$ be a variation of $\ZZ$-mixed Hodge structure on a Zariski-open subset $S$
of a complex manifold $\Sb$. We call $\shV$ \define{admissible with respect to $\Sb$} if
$\shV \tensor \QQ$ is admissible in the sense of Kashiwara \cite{Kashiwara} (where
admissibility is defined by a curve test). It is clear from this definition that
admissibility is preserved under holomorphic maps $f \colon \Sb' \to \Sb$ with the
property that $f^{-1}(S)$ is dense in $\Sb'$. Moreover, duals and tensor products of admissible
variations of mixed Hodge structure are again admissible; this is proved in the
appendix to \cite{SZ}.

By work of Saito \cite{Saito-ANF}, admissibility can also be phrased in terms of
mixed Hodge modules: $\shV \tensor \QQ$ defines a mixed Hodge module on $S$, and is
admissible if and only if that mixed Hodge module can be extended to $\Sb$.

\section{The locus of Hodge classes}
\label{s.LocusOfHodge}

We now turn to the proof of Theorem~\ref{thm:main}. Throughout, we let $\shV$ be a
variation of mixed Hodge structure over $S$, admissible with respect to $\Sb$. We
assume that $\shV$ is graded polarized, and that the local systems $W_m \shV$ making
up the weight filtration are defined over $\ZZ$, with $\Gr_m^W \shV$ torsion free.

To begin with, we can replace $\shV$ by $W_0 \shV$, and assume without loss of
generality that $\shV$ is of weight $\leq 0$. We then have 
\[
	\Hdg(\shV) = T(\shVZ) \cap E(F^0 \shVO).
\]

The next step is to prove a more general version of Theorem~\ref{thm:BP}. Recall that
a \define{generalized normal function} $\nu$ is an extension, in the category of
variations of mixed Hodge structure, of $\ZZ(0)$ by a variation of mixed Hodge
structure $\shH$, all of whose weights are $\leq -1$.  It is said to be
\define{admissible} if the corresponding variation is admissible. At each point $s
\in S$, the extension determines a point $\nu(s) \in \Ext_{\MHS}^1 \bigl( \ZZ(0), H_s
\bigr)$; the \define{zero locus} $Z(\nu)$ of the generalized normal function is by
definition the set of points where $\nu(s) = 0$. We let
\[
	\NF(S, \shH) = \Ext_{\VMHS(S)}^1 \bigl( \ZZ(0), \shH \bigr)
\]
denote the group of generalized normal functions.

\begin{prop} \label{prop:BPgen}
Let $\nu$ be an admissible generalized normal function on $S$.  Then the
closure of $Z(\nu)$ in $\Sb$ is an analytic subset.
\end{prop}

\begin{proof}
Let $\shV$ be the corresponding admissible variation of mixed Hodge structure, and
$\shH = W_{-1} \shV$. If $\shH$ is pure, then the result follows from
Theorem~\ref{thm:BP}. Otherwise, we let $m \leq -1$ be the smallest integer for which $\Gr_m^W \shV
\neq 0$. Define $\shV' = \shV / W_m \shV$, and let $\nu_0$ be the corresponding
generalized normal function induced on $\shV'$ by $\nu$. Note that we have $Z(\nu) \subseteq Z(\nu_0)$.

Let $S_0$ denote the regular locus of an irreducible component of $Z(\nu_0)$. By induction, we know that the
closure of $S_0$ inside of $\Sb$ is analytic; let $\pi \colon \Sb_0 \to \Sb$ be a
resolution of singularities of the closure that is an isomorphism over $S_0$. Since
$\pi$ is proper, we are allowed to replace $\Sb$ by $\Sb_0$ and $\nu$ by its pullback
to $S_0$; we may therefore assume from the beginning that $\nu_0 = 0$. Now the
exact sequence
\begin{diagram}
0 &\rTo& \NF \bigl( S, W_m \shH \bigr) &\rTo& \NF(S, \shH) &\rTo&
	\NF \bigl( S, \shH / W_m \shH \bigr)
\end{diagram}
shows that $\nu$ induces a generalized normal function $\nu' \in \NF \bigl( S, W_m
\shH \bigr)$. Since $W_m \shH$ is pure of weight $m$, we conclude from
Theorem~\ref{thm:BP} that $Z(\nu')$ has an analytic closure inside $\Sb$; but clearly
$Z(\nu) = Z(\nu')$, and so the assertion follows.
\end{proof}

We are now ready to prove Theorem~\ref{thm:main} in general.
\begin{proof}[Proof of Theorem~\ref{thm:main}]
Let $\shV$ be the admissible variation of mixed Hodge structure; as explained above,
we may suppose that it has weights $\leq 0$. For any point $s \in S$, let $V_s$ be
the corresponding mixed Hodge structure; then we have an exact sequence
\begin{equation} \label{eq:Hdg}
\begin{diagram}
0 &\rTo& \Hdg(V_s) &\rTo& \Hdg \bigl( \Gr_0^W V_s \bigr) &\rTo&
	\Ext_{\MHS}^1 \bigl( \ZZ(0), W_{-1} V_s \bigr).
\end{diagram}
\end{equation}
It follows that the locus of Hodge classes for $\shV$ is embedded into that for
$\Gr_0^W \shV$. Let $Z$ be an irreducible component of $\Hdg(\shV)$, and let $Y$ be
the irreducible component of $\Hdg \bigl( \Gr_0^W \shV \bigr)$ containing $Z$.
By the theorem of Cattani, Deligne, and Kaplan \cite{CDK}, $Y$ can be extended to an analytic space $\Yb$
that is proper and finite over $\Sb$. Let $\mu \colon \Yb' \to \Yb$ be a resolution
of singularities of the analytic space $\Yb$
and denote
by $\shV'$ the pullback of $\shV$ to $Y$. 

By construction, we have a section $\ZZ(0) \to \Gr_0^W \shV'$. It induces
a generalized normal function $\nu' \in \NF(Y, \shH')$, where $\shH' = W_{-1} \shV'$.
Moreover, it is clear from \eqref{eq:Hdg} that $Z = Z(\nu')$. Since $\nu'$
is easily seen to be admissible with respect to $\Yb'$, we conclude from
Proposition~\ref{prop:BPgen} that the closure of $Z(\nu')$ in $\Yb'$ is analytic.
Because $\mu$ is proper, it follows that $Z$ has an analytic closure inside of $\Yb$;
this completes the proof.
\end{proof}

\section{The split locus}

The proof of Theorem~\ref{thm:split} is similar to that of Theorem~\ref{thm:main}. 
\begin{proof}
  It suffices to prove the statement with coefficients in $\QQ$. So
  let $\shV$ be an admissible variation of mixed Hodge structure on
  $S$, where $S$ is Zariski-open in a complex manifold $\Sb$. Let $m$
  be the largest integer for which $\Gr_m^W \shV \neq 0$. By
  induction, we know that the locus of points $s \in S$ where $W_{m-1}
  V_s$ splits over $\QQ$ has an analytic closure inside of
  $\Sb$. Arguing as before, we may therefore assume from the beginning
  that $W_{m-1} \shV$ is split. Now $\shV$ determines an element of
\begin{align*}
	\Ext_{\VMHS(S)}^1 \bigl( \Gr_m^W \shV, W_{m-1} \shV \bigr) 
		&\simeq \bigoplus_{k < m} 
		\Ext_{\VMHS(S)}^1 \bigl( \Gr_m^W \shV, \Gr_k^W \shV \bigr) \\
		&\simeq \bigoplus_{k < m}
		\Ext_{\VMHS(S)}^1 \bigl( \QQ(0), (\Gr_m^W \shV)^{\vee} \tensor \Gr_k^W \shV \bigr).
\end{align*}
Since admissibility is preserved under tensor products, the problem is reduced to the
case of admissible higher normal functions; applying Theorem~\ref{thm:BP} completes
the proof.
\end{proof}



\end{document}